\documentclass{article}
\usepackage{amsmath}
\usepackage{amsfonts}
\usepackage{amssymb}

\newtheorem{theorem}{Theorem}

\newtheorem{corollary}[theorem]{Corollary}

\newtheorem{lemma}[theorem]{Lemma}

\newtheorem{proposition}[theorem]{Proposition}

\input{tcilatex}
\begin{document}

\begin{center}
{\LARGE An Integral Kernel }

{\LARGE for Weakly Pseudoconvex Domains}\footnote{%
}\footnote{%
Part of this research was begun while the author was visiting at the
University of Utah in Winter 2011. The author would like to thank the
Mathematics Department of that University for the hospitality and support
extended to him.}

\bigskip

R. Michael Range\footnote{%
Department of Mathematics, State University of New York at Albany, Albany,
NY 12222. $e-mail:range@math.albany.edu$}

\vspace{0.5in}

\textbf{ABSTRACT\bigskip }
\end{center}

\textit{A new explicit construction of Cauchy-Fantappi\'{e} kernels is
introduced for an arbitrary weakly pseudoconvex domain with smooth boundary.
While not holomorphic in the parameter, the new kernel reflects the complex
geometry and the Levi form of the boundary. Some estimates are obtained for
the corresponding integral operator, which provide evidence that this kernel
and related constructions give useful new tools for complex analysis on this
general class of domains. }

\begin{center}
\bigskip
\end{center}

\section{Introduction}

The well-known Bochner-Martinelli kernel $K_{BM}(\zeta ,z)$ is a natural
generalization of the familiar Cauchy kernel to higher dimensions. It leads
to a corresponding integral representation formula for holomorphic functions
on \emph{arbitrary} smoothly bounded domains in $\mathbb{C}^{n}$. However,
it lacks critical properties, and this limits its applicability. For once,
it is not holomorphic in the parameter $z$ when $n\geq 2.$ Furthermore, its
singularity is isotropic, and thus it does not at all reflect the
non-isotropic \emph{complex }geometry of boundaries of domains in higher
dimension.

In the late 1960s G. Henkin and E. Ramirez, independently, introduced new
integral representation formulas on \emph{strictly }pseudoconvex domains
which overcame the shortcomings of the Bochner-Martinelli kernel. These new
tools rapidly led to proofs of numerous results concerning the boundary
behavior of holomorphic functions and related objects on strictly
pseudoconvex domains. (See Ra86 for a systematic exposition.) In particular,
they allowed to prove \emph{pointwise} estimates for solutions of the
Cauchy-Riemann equations, such as estimates in supremum norm and in H\"{o}%
lder norms, which were not accessible by the classical $L^{2}$ - methods of
J. J. Kohn, L. H\"{o}rmander, and others. Note that in dimension one every
smoothly bounded domain is trivially strictly pseudoconvex. Such domains
therefore provide a natural setting for generalizing results known in
dimension one to higher dimensions.

Of course, when $n\geq 2$, not every domain of holomorphy with smooth
boundary is \emph{strictly} pseudoconvex. Such more general domains are just
(weakly) pseudoconvex, i.e., the Levi form associated to the boundary is
only positive \emph{semi-definite} rather than positive definite as in the 
\emph{strict} case. This more general case has been investigated for well
over 40 years, and it continues to present major challenges. As conjectured
by J. J. Kohn, the right notion of \emph{finite type} (see Ko72, Ko 79,
Da82) turned out to be central in the $L^{2}$ theory of the $\overline{%
\partial }-$Neumann problem, providing necessary and sufficient conditions
for the existence of \emph{subelliptic }estimates (Ca87). However, attempts
to generalize kernel methods to this setting and obtain, for example,
pointwise estimates for solutions of the $\overline{\partial }$-equation,
have had only limited success (see Ra78, Ra90, Cu97, DF99, DFF99). In fact,
there is a fundamental obstruction to extending these methods to this
setting, as follows. The Henkin/Ramirez construction made essential use of
an explicit \emph{holomorphic} \emph{support function}, which in the
strictly pseudoconvex case is given locally by the quadratic Levi
polynomial. Holomorphic support functions exist also on Euclidean convex
domains, but as shown by an example of Kohn and L. Nirenberg [KoNi72], they
do not exist in general for pseudoconvex domains of finite type.
Consequently the integral kernel methods seemed to have reached their limit,
and a more complete understanding of pointwise estimates in the theory of
the $\overline{\partial }-$Neumann problem in arbitrary dimensions has
remained elusive for quite some time.

In this paper we introduce a \emph{non-holomorphic} local modification $\Phi
(\zeta ,z)$ of the Levi polynomial for an arbitrary weakly pseudoconvex
domain $D\subset \subset \mathbb{C}^{n},$ and a corresponding global kernel
generating form \medskip

\begin{equation*}
W^{S}(\zeta ,z)=\frac{\sum_{j=1}^{n}s_{j}(\zeta ,z)d\zeta _{j}}{S(\zeta ,z)}%
\text{ on }bD\times D,
\end{equation*}%
with $S(\zeta ,z)=\Phi (\zeta ,z)\ $for $z$ close to $\zeta ,$ and where the
coefficients $s_{j}$ satisfy $S(\zeta ,z)=$\emph{\ }$\sum_{j=1}^{n}s_{j}(%
\zeta _{j}-z_{j})$. While $\Phi $ and $W^{S}$ are not holomorphic in $z$,
they satisfy critical basic properties which open the door to significant
applications and new results. In particular, we shall prove the following
properties.

a)\emph{\ The form }$\overline{\partial _{z}}\Phi (\zeta ,z)$\emph{\ has a
zero at }$z=\zeta $\emph{\ whose order is carefully controlled.}

b) $\Phi $\emph{\ satisfies precise uniform estimates from below somewhat
weaker than those familiar in the strictly pseudoconvex case, and which
involve explicitly the eigenvalues of the Levi form of a defining function }$%
r$\emph{\ for }$D$\emph{.}

c) \emph{In case the domain is strictly pseudoconvex, }$\Phi $ \emph{and }$%
W^{S}$\emph{---while not holomor- phic---satisfy the classical estimations
known in that case.}

d) $\Phi $\emph{\ satisfies the same symmetry properties that have been
successfully used on strictly pseudoconvex domains in earlier work.\medskip }

The resulting Cauchy-Fantappi\'{e} kernel $\Omega _{0}(W^{S})=(2\pi
i)^{-n}W^{S}\wedge (\overline{\partial }_{\zeta }W^{S})^{n-1}$ yields a
Cauchy-type integral formula for holomorphic functions on a weakly
pseudoconvex domain $D$ which reflects the complex geometry of the boundary.
In particular, it has a singularity of order one in the complex normal
direction (just like the standard one dimensional Cauchy kernel), and a
singularity of order two in each of the complex tangential directions.
Furthermore, it will be shown that $\overline{\partial _{z}}\Omega
_{0}(W^{S})$ is in some sense more regular than the corresponding form for
the Bochner-Martinelli kernel. These features suggest that $\Omega
_{0}(W^{S})$ might be a useful new tool for complex analysis on weakly
pseudoconvex domains. For example, the author has used $\Omega _{0}(W^{S})$
and its higher order versions $\Omega _{q}(W^{S}),$ $0\leq q\leq n$, in the
theory of the $\overline{\partial }-$Neumann problem. A major new result he
obtained is a \emph{pointwise} analogon of the classical basic $L^{2}$
estimate of Morrey and Kohn (see FoKo72). Furthermore, some additional
results suggest that it might be possible to develop a suitable version of
Kohn's theory of subelliptic multipliers (see Ko79, Siu10) in the integral
representation setting involving the kernels $\Omega _{q}(W^{S})$ and their
variants. Such techniques might then allow to prove suitable H\"{o}lder
estimates on pseudoconvex domains of finite type. A preliminary report
discussing these applications has been published on arXiv (Ra11).

\section{Local construction of the support function}

We assume that $D$ is a bounded pseudoconvex domain in $\mathbb{C}^{n}$ with 
$C^{k}$ boundary $bD$ ( $k$ $\geq 3$), and we choose a $C^{k}$ defining
function $\varphi $ for $bD$ defined on a neighborhood $U=U(bD)$ of $bD$. In
general the level surfaces $M_{-\delta }=\{z:\varphi (z)=-\delta \}$ will
not be Levi pseudoconvex for $\delta >0$.

\begin{proposition}
There exists $C>0$ and $U$, such that for all $\zeta \in \overline{D}\cap U$
the Levi form $\mathcal{L}$ of the defining function $r(z)=\varphi (z)\exp
(-C\left\vert z\right\vert ^{2})$ satisfies 
\begin{equation*}
\mathcal{L}(r,\zeta ;t)=\sum_{j,k=1}^{n}\frac{\partial ^{2}r}{\partial \zeta
_{j}\partial \overline{\zeta _{k}}}(\zeta )t_{j}\overline{t_{k}}\geq 0\text{
for all }t\in \mathbb{C}^{n}\text{ with}\sum_{j=1}^{n}\frac{\partial r}{%
\partial \zeta _{j}}(\zeta )t_{j}=0.
\end{equation*}
\end{proposition}

As the proof will show, the level surfaces $M_{r(\zeta )}$ of $r$ are
actually strictly pseudoconvex for $\zeta \in D\cap U$, but the resulting
estimates are not uniform in $\zeta $ as $r(\zeta )\rightarrow 0$, unless $%
bD $ is strictly pseudoconvex to begin with.

\textit{Proof.} By Theorem 2 in [Ra81], there exist $C>0$ and $0<\eta <1$
such that $\rho =-(-\varphi \exp (-C\left\vert z\right\vert ^{2}))^{\eta }$
is strictly plurisubharmonic on $U\cap D$ for $U$ sufficiently small. In
particular, the level surfaces of $\rho $ close to $bD$ are \emph{strictly
pseudoconvex}. Let $r=$ $\varphi \exp (-C\left\vert z\right\vert ^{2})$.
Since $\rho (z)=-\delta ^{\eta }$ if and only if $r(z)=-\delta $, it follows
that the level surfaces $M_{-\delta }$ of $r$ are (strictly) pseudoconvex as
well for $\delta $ sufficiently small, and the desired result follows.

We now fix this particular global defining function $r$. After shrinking $U$%
, we may assume that for a fixed $k\geq 3$ the function $r$ has a bounded $%
C^{k}$ norm $\left\vert r\right\vert _{k}$ over $U$.

For $\zeta \in U(bD)$ let 
\begin{equation*}
F^{(r)}(\zeta ,z)=\sum_{j}\frac{\partial r}{\partial \zeta _{j}}(\zeta
)(\zeta _{j}-z_{j})-\frac{1}{2}\sum_{j,k}\frac{\partial ^{2}r}{\partial
\zeta _{j}\partial \zeta _{k}}(\zeta )(\zeta _{j}-z_{j})(\zeta _{k}-z_{k})
\end{equation*}%
be the usual Levi polynomial of $r,$ which is a quadratic holomorphic
polynomial in $z$. (See Ra86) The following well known equation is a direct
consequence of the 2nd order Taylor expansion of $r(z)$ at $\zeta $. 
\begin{equation*}
2\func{Re}\left[ F^{(r)}(\zeta ,z)-r(\zeta )\right] =-r(\zeta )-r(z)+%
\mathcal{L}(r,\zeta ;\zeta -z)+O(\left\vert \zeta -z\right\vert ^{3}).
\end{equation*}%
We now define

\begin{equation*}
\Phi _{K}(\zeta ,z)=F^{(r)}(\zeta ,z)-r(\zeta )+K\left\vert \zeta
-z\right\vert ^{3}\text{ ,}
\end{equation*}%
where $K>0$ is a large constant to be suitably chosen later on.

While $\left\vert \zeta -z\right\vert ^{2}$ is smooth in $(\zeta ,z)$, the
term $\left\vert \zeta -z\right\vert ^{3}$ which appears in $\Phi _{K}$ is
of class $C^{2}$ in general and smooth only at points $(\zeta ,z)$ with $%
\zeta \neq z.$ For $j=0,1,2,...$ we denote by $\mathcal{E}_{j}=\mathcal{E}%
_{j}(\zeta ,z)$ a smooth function or form which satisfies $\left\vert 
\mathcal{E}_{j}\right\vert \leq const.\left\vert \zeta -z\right\vert ^{j}$
for a constant that is independent of $\zeta $ and $z$; similarly, we use $%
\mathcal{E}_{j}^{\#}$ to denote \emph{bounded} functions or forms which are
smooth for $\zeta ,z$ with $\zeta \neq z$, and which satisfy an estimate $%
\left\vert \mathcal{E}_{j}^{\#}\right\vert \leq const.\left\vert \zeta
-z\right\vert ^{j}$. Then $\left\vert \zeta -z\right\vert ^{3}=\mathcal{E}%
_{3}^{\#}$, and one readily verifies that if $D^{l}$ is a partial derivative
of order $l$ with respect to $\zeta $ and/or $z$, then $D^{l}\left\vert
\zeta -z\right\vert ^{3}=\mathcal{E}_{3-l}^{\#}$ for $l=1,2,3.$ We point out
that the precise expressions of both $\mathcal{E}_{j}$ and $\mathcal{E}%
_{j}^{\#}$ may differ from formula to formula, and even within the same
formula. The emphasis is on keeping track of the order of vanishing of
terms, not of their exact expressions. This will become most relevant in
section 5.

Since $F^{(r)}-r(\zeta )$ is holomorphic in $z$, it follows that $\overline{%
\partial _{z}}\Phi _{K}(\zeta ,z)=\mathcal{E}_{2}^{\#}$. For the other
derivatives, fix $P\in bD$ and introduce a $C^{k-1}$ orthonormal frame $%
\{\omega _{1},$ $\omega _{2},...,\omega _{n}\}$ for $(1,0)$ forms on a
sufficiently small neighborhood $V(P),$ with $\omega _{n}=\nu \partial r$
for some function $\nu (\zeta )>0$ on $V$, so that $\left\Vert \omega
_{n}\right\Vert =1$ on $V$. Let $\{L_{1},...,L_{n}\}$ be the corresponding
dual frame for $(1,0)$ vector fields on $V.$ Then $L_{j}(r)=0$ for $j<n$,
and $L_{n}=\gamma (\zeta )\sum_{k}\overline{\partial r/\partial \zeta _{k}}%
~\partial /\partial \zeta _{k}$ for some $\gamma (\zeta )>0.$ The vector
fields $L_{j}$ act in $\zeta ;$ we use the notation $L_{j,z}$ and $\overline{%
L_{j,z}}$ if differentiation is taken with respect to $z.$

\begin{proposition}
The following estimates hold for the derivatives of $\Phi _{K}$:
\end{proposition}

i) $\overline{L_{j,z}}\Phi _{K}(\zeta ,z)$ $=\mathcal{E}_{2}^{\#}$ \textit{%
for} $j=1,...,n$;

ii) $L_{j,z}\Phi _{K}(\zeta ,z)$ $=\mathcal{E}_{1}^{\#}$ \textit{for} $j<n$;

iii) $L_{n,z}\Phi _{K}(\zeta ,z)\neq 0.$

\textit{Proof.} We already noted i). For ii), note that if $%
L_{j}=\sum_{k}a_{jk}(\zeta )~\partial /\partial \zeta _{k},$ then $%
L_{j,z}(F^{(r)}-r(\zeta ))=-\sum_{k}a_{jk}(\zeta )~\partial r/\partial \zeta
_{k}+\mathcal{E}_{1}=\mathcal{E}_{1}$ for $j<n.$ This implies ii). Finally,
since $L_{n,z}F^{(r)}(\zeta ,z)=-\gamma (\zeta )\sum_{k}\left\vert \partial
r/\partial \zeta _{k}\right\vert ^{2}+\mathcal{E}_{1}$, iii) follows as well.

One also has the following approximate symmetry, which follows directly from
the known result in case $K=0$ (see Ra 86).

\begin{equation*}
\Phi _{K}(\zeta ,z)-\overline{\Phi _{K}(z,\zeta )}=\mathcal{E}_{3}.\medskip
\end{equation*}

\section{Estimations for the support function}

Next we prove that $\left\vert \Phi _{K}(\zeta ,z)\right\vert $ is precisely
controlled from below, as follows. We use the convention that $A\succsim B$
means that there exists a constant $c>0$, so that $\left\vert A\right\vert
\geq c\left\vert B\right\vert $ for all the values of the relevant variables
under consideration.

For $\zeta \in U$ consider the level surface $M_{r(\zeta )}$ of $r$ through
the point $\zeta $. We introduce the orthogonal projection $\pi _{\zeta
}^{t}:\mathbb{C}^{n}\rightarrow T_{\zeta }^{1,0}(M_{r(\zeta )})\subset 
\mathbb{C}^{n}$, where $T_{\zeta }^{1,0}(\mathbb{C}^{n})$ is identified with 
$\mathbb{C}^{n}$ via the standard basis $\{\frac{\partial }{\partial \zeta
_{1}},...,\frac{\partial }{\partial \zeta _{n}}\}$.

\begin{theorem}
The neighborhood $U$, the constant $K$, and $\varepsilon >0$ can be chosen
so that for all $\zeta ,z$ $\in \overline{D}\cap U$ with $\left\vert \zeta
-z\right\vert <\varepsilon $ one has
\end{theorem}

\begin{eqnarray*}
\left\vert \Phi _{K}(\zeta ,z)\right\vert &\gtrsim &\left[ \left\vert \func{%
Im}F^{(r)}(\zeta ,z)\right\vert +\left\vert r(\zeta )\right\vert +\left\vert
r(z)\right\vert +\right. \\
&&\left. +\mathcal{L}(r,\zeta ;\pi _{\zeta }^{t}(\zeta -z))+K\left\vert
\zeta -z\right\vert ^{3}\right]
\end{eqnarray*}%
\medskip Note that by pseudoconvexity and by the special choice of the
defining function $r$ one has 
\begin{equation*}
\mathcal{L}(r,\zeta ;\pi _{\zeta }^{t}(\zeta -z))\geq 0\text{ for all }\zeta
\in \overline{D}\cap U\text{,}
\end{equation*}%
so all terms on the right side in the estimation of $\left\vert \Phi
_{K}(\zeta ,z)\right\vert $ are nonnegative!

As in the familiar strictly pseudoconvex case, $r(\zeta )$ and $\func{Im}%
F^{(r)}(\zeta ,z)$ can be used as coordinates in a $C^{k-2}$ real coordinate
system in a neighborhood $B(z,\delta )$ of a fixed point $z\in U$ provided $%
\delta >0$ is sufficiently small. The crux of the estimate in the Theorem is
that $\Phi _{K}$ is of order $1$ in the complex \emph{normal} direction,
while the Levi form completely controls $\Phi _{K}$ from below in the
complex \emph{tangential} directions.

\textit{Proof. }Decompose $\zeta -z=\pi _{\zeta }^{t}(\zeta -z)+\pi _{\zeta
}^{n}(\zeta -z)$, and notice that $\left\vert \pi _{\zeta }^{n}(\zeta
-z)\right\vert =O(\left\vert <\partial r(\zeta ),\zeta -z>\right\vert )$,
where $<\partial r(\zeta ),\zeta -z>=$ $\sum_{j}\frac{\partial r}{\partial
\zeta _{j}}(\zeta )(\zeta _{j}-z_{j})$.

Since%
\begin{equation*}
<\partial r(\zeta ),\zeta -z>=\Phi _{K}+r(\zeta )+\mathcal{E}%
_{2}-K\left\vert \zeta -z\right\vert ^{3}\text{,}
\end{equation*}%
it follows that 
\begin{eqnarray*}
\mathcal{L}(r,\zeta ;\zeta -z) &=&\mathcal{L}(r,\zeta ;\pi _{\zeta
}^{t}(\zeta -z))+\mathcal{E}_{1}\left\vert <\partial r(\zeta ),\zeta
-z>\right\vert \geq \\
&\geq &\mathcal{L}(r,\zeta ;\pi _{\zeta }^{t}(\zeta -z))-\mathcal{E}%
_{1}[\left\vert \Phi _{K}\right\vert +\left\vert r(\zeta )\right\vert
+K\left\vert \zeta -z\right\vert ^{3}]-\mathcal{E}_{3}\text{.}
\end{eqnarray*}%
We choose $\varepsilon >0$ so small that the above $\mathcal{E}_{1}$ term
satisfies $\left\vert \mathcal{E}_{1}\right\vert <1/2$ for $\left\vert \zeta
-z\right\vert <\varepsilon .$ It follows that for $\zeta ,z\in \overline{D}%
\cap U$ with $\left\vert \zeta -z\right\vert <\varepsilon $ one has 
\begin{eqnarray*}
2\func{Re}\Phi _{K} &=&2\func{Re}\left[ F^{(r)}(\zeta ,z)-r(\zeta )\right]
+2K\left\vert \zeta -z\right\vert ^{3}\geq \\
&\geq &\left\vert r(\zeta )\right\vert +\left\vert r(z)\right\vert +\mathcal{%
L}(r,\zeta ;\pi _{\zeta }^{t}(\zeta -z))+\frac{3}{2}K\left\vert \zeta
-z\right\vert ^{3}+ \\
&&-1/2(\left\vert \Phi _{K}\right\vert +\left\vert r(\zeta )\right\vert
)-A\left\vert \zeta -z\right\vert ^{3}
\end{eqnarray*}%
for a certain constant $A>0.$ Now choose $K\geq 2A$; it then follows that 
\begin{eqnarray*}
2\func{Re}\Phi _{K} &\geq &\left\vert r(\zeta )\right\vert +\left\vert
r(z)\right\vert +\mathcal{L}(r,\zeta ;\pi _{\zeta }^{t}(\zeta
-z))+K\left\vert \zeta -z\right\vert ^{3}+ \\
&&-1/2(\left\vert \Phi _{K}\right\vert +\left\vert r(\zeta )\right\vert ).
\end{eqnarray*}%
After rearranging, it follows that for $\zeta ,z\in \overline{D}\cap U$ with 
$\left\vert \zeta -z\right\vert <\varepsilon $ one has 
\begin{equation*}
2\func{Re}\Phi _{K}+\left\vert \Phi _{K}\right\vert /2\geq \left\vert
r(\zeta )\right\vert /2+\left\vert r(z)\right\vert +\mathcal{L}(r,\zeta ;\pi
_{\zeta }^{t}(\zeta -z))+K\left\vert \zeta -z\right\vert ^{3}.
\end{equation*}

Since $\left\vert \Phi _{K}\right\vert \succeq 2\left\vert \func{Re}\Phi
_{K}\right\vert +2\left\vert \func{Im}\Phi _{K}\right\vert \geq 2\func{Re}%
\Phi _{K}+2\left\vert \func{Im}F^{r}\right\vert $, the preceding inequality
readily gives the estimate stated in the Theorem.$\blacksquare $

\textbf{Remark. }If $bD$ is strictly pseudoconvex, then 
\begin{equation*}
L(r,\zeta ;\pi _{\zeta }^{t}(\zeta -z))\geq c\left\vert \pi _{\zeta
}^{t}(\zeta -z)\right\vert ^{2}\text{ for some }c>0.
\end{equation*}%
Note that 
\begin{equation*}
\left\vert \zeta -z\right\vert ^{2}\leq \left\vert \pi _{\zeta }^{t}(\zeta
-z)\right\vert ^{2}+\mathcal{E}_{1}\left\vert <\partial r(\zeta ),\zeta
-z>\right\vert
\end{equation*}%
By estimating $\left\vert <\partial r(\zeta ),\zeta -z>\right\vert $ as
before, one obtains 
\begin{equation*}
\left\vert \pi _{\zeta }^{t}(\zeta -z)\right\vert ^{2}\geq \left\vert \zeta
-z\right\vert ^{2}-A_{1}\left\vert \zeta -z\right\vert [\left\vert \Phi
_{K}\right\vert +\left\vert r(\zeta )\right\vert )+K\left\vert \zeta
-z\right\vert ^{3}]-A_{2}\left\vert \zeta -z\right\vert ^{3},
\end{equation*}%
where $A_{1},$ $A_{2}$ are certain positive constants. Given $K\geq 0$, one
can then choose $\varepsilon =\varepsilon (K)$ so small that 
\begin{equation*}
\left\vert \pi _{\zeta }^{t}(\zeta -z)\right\vert ^{2}\succeq \left\vert
\zeta -z\right\vert ^{2}-1/2(\left\vert \Phi _{K}\right\vert +\left\vert
r(\zeta )\right\vert )\text{ for }\left\vert \zeta -z\right\vert
<\varepsilon .
\end{equation*}%
This readily implies the stronger estimate

\begin{equation*}
\left\vert \Phi _{K}(\zeta ,z)\right\vert \succsim \left\vert \func{Im}%
F^{(r)}(\zeta ,z)\right\vert +\left\vert r(\zeta )\right\vert +\left\vert
r(z)\right\vert +\widetilde{c}\left\vert \zeta -z\right\vert ^{2}
\end{equation*}%
for any $K\geq 0$ and for $\zeta ,z\in \overline{D}\cap U$ with $\left\vert
\zeta -z\right\vert \leq \varepsilon $, provided $\varepsilon $ is chosen
sufficiently small.\medskip\ This includes the familiar estimate for the
Levi polynomial in the strictly pseudoconvex case, which in the literature
has usually been obtained for a \emph{strictly plurisubharmonic} defining
function. The proof here gives the result for an arbitrary defining function.

Most importantly, the\textit{\ third }order correction term allows to prove
the following more delicate estimate which is critical for the estimations
of Cauchy-Fantappi\'{e} kernels involving the support function $\Phi _{K}$.

\begin{proposition}
Fix $z\subset \overline{D}\cap U$ and let $\lambda _{1}(z),...,\lambda
_{n-1}(z)$ be the eigenvalues of the Levi form of the chosen defining
function $r$ at the point $z.$ There exists a unitary change of coordinates
in the $\zeta $ variables (in dependence of $z$), so that in the new
coordinates one has the estimate
\end{proposition}

\begin{equation*}
\left\vert \Phi _{K}(\zeta ,z)\right\vert \succsim \left\vert \func{Im}%
F^{(r)}(\zeta ,z)\right\vert +\left\vert r(\zeta )\right\vert +\left\vert
r(z)\right\vert +\sum_{j=1}^{n-1}\lambda _{j}(z)\left\vert \zeta
_{j}-z_{j}\right\vert ^{2}+K/2\left\vert \zeta -z\right\vert ^{3}~
\end{equation*}%
\textit{for }$\zeta ,z\in \overline{D}\cap U$\textit{\ with }$\left\vert
\zeta -z\right\vert <\delta $, \textit{provided }$K$ \textit{is sufficiently
large and }$\delta $\textit{\ is sufficiently small}.\textit{\ All relevant
constants can be chosen to be independent of }$\zeta $ and $z\in \overline{D}%
\cap U.$

\textit{Proof.} Fix $z$, and choose the orthonormal frame $L_{1},...,L_{n}$
so that $L_{1},...,L_{n-1}$ form an orthonormal basis for $T_{\zeta
}^{1,0}(M_{r(\zeta )})$ for $\left\vert \zeta -z\right\vert <\delta \leq
\varepsilon $ which diagonalizes the Levi form restricted to $%
T_{z}^{1,0}(M_{r(z)})$ at the point $\zeta =z$. Note that this is a
condition at the single point $\zeta =z$; in general, there is no smooth
frame which diagonalizes the Levi form in a neighborhood of $z$. By
pseudoconvexity and the choice of the defining function the eigenvalues $%
\lambda _{j}(z),$ $j=1,...,n-1$, are nonnegative. After a unitary change of
coordinates in $\zeta _{1},...,\zeta _{n},$ one can assume that $\left.
L_{j}\right\vert _{z}=\sqrt{2}\left. \frac{\partial }{\partial \zeta _{j}}%
\right\vert _{z}$, and hence $\left. L_{j}\right\vert _{\zeta }=\sqrt{2}%
\left. \frac{\partial }{\partial \zeta _{j}}\right\vert _{\zeta }+V_{j}$,
where the coefficients of $V_{j}$ are of type $\mathcal{E}_{1}$. We call
such coordinates $z-$ \emph{diagonalizing}. It then follows that with
respect to these particular coordinates one has 
\begin{equation*}
\mathcal{L}(r,z;\pi _{z}^{t}(\zeta -z))=\sum_{j=1}^{n-1}\lambda
_{j}(z)\left\vert \zeta _{j}-z_{j}\right\vert ^{2}.\text{ }
\end{equation*}%
Since the coefficients of the Levi form are smooth in $\zeta $, one obtains%
\begin{equation*}
\mathcal{L}(r,\zeta ;\pi _{\zeta }^{t}(\zeta -z))=\sum_{j=1}^{n-1}\lambda
_{j}(z)\left\vert \zeta _{j}-z_{j}\right\vert ^{2}+\mathcal{R}(\zeta ,z)%
\text{,}
\end{equation*}%
where the error term $\mathcal{R}(\zeta ,z)$ is of type $\mathcal{E}_{3}.$
Now choose $K$ in the definition of $\Phi _{K}$ so large that this error
term satisfies $\left\vert \mathcal{R}(\zeta ,z)\right\vert \leq \frac{K}{2}%
\left\vert \zeta -z\right\vert ^{3}$. The desired estimate then follows from
the estimate in Theorem 3.$\blacksquare $

There is a related estimate involving the corresponding dual frame, as
follows.

\begin{proposition}
With $z\subset \overline{D}\cap U$ fixed, and the frame $\{L_{1},...,L_{n}\}$
and coordinates chosen as in the proof of the preceding proposition, one has 
\begin{equation*}
\partial r\wedge \overline{\partial }r\wedge \partial \overline{\partial }%
r(\zeta )=\gamma (\zeta )\omega _{n}\wedge \overline{\omega }_{n}\wedge %
\left[ 1/2\sum_{j=1}^{n-1}\lambda _{j}(z)d\zeta _{j}\wedge \overline{d\zeta
_{j}}+\Omega _{1}\right] ,
\end{equation*}%
where $\gamma (\zeta )\neq 0$ and the $(1,1)$ form $\Omega _{1}$ has smooth $%
\mathcal{E}_{1}$ coefficients.
\end{proposition}

\textit{Proof. }Let $\{\omega _{1},...,\omega _{n}\}$ be the corresponding
dual frame for $(1,0)$ forms on $\{\zeta :\left\vert \zeta -z\right\vert
<\delta \}$. By the special choice of the frame one has\textit{\ }$\partial 
\overline{\partial }r(z)=\sum_{j=1}^{n-1}\lambda _{j}(z)\omega _{j,z}\wedge 
\overline{\omega _{j,z}}+\Lambda _{1}\wedge \omega _{n}+\Lambda _{2}\wedge 
\overline{\omega _{n}}$, with certain $1$-forms $\Lambda _{1}$ and $\Lambda
_{2}$. Note that $\omega _{j,z}=1/\sqrt{2}\left. d\zeta _{j}\right\vert _{z}$%
. Since $\partial \overline{\partial }r(\zeta )$ and the frame\emph{\ }are
smooth, letting the base point $\zeta =z$ in the forms vary while freezing
the eigenvalues at the point $z$, it follows that 
\begin{equation*}
\partial \overline{\partial }r(\zeta )=\left[ 1/2\sum_{j=1}^{n-1}\lambda
_{j}(z)d\zeta _{j}\wedge \overline{d\zeta _{j}}\right] +\Omega _{1}+\Lambda
_{1}\wedge \omega _{n}+\Lambda _{2}\wedge \overline{\omega _{n}},
\end{equation*}%
for a $(1,1)$-form $\Omega _{1}$ with $\mathcal{E}_{1}$ coefficients. The
proposition then follows by wedging the last equation with $\partial r\wedge 
\overline{\partial }r=\gamma (\zeta )~\omega _{n}\wedge \overline{\omega }%
_{n}$. $\blacksquare $

\begin{corollary}
With $z\in \overline{D}\cap U$ fixed and $\zeta $ the corresponding $z$%
-diagonalizing coordinates as above, one has the following representation
for $\zeta \in \overline{D}\cap U$ and $0<\left\vert \zeta -z\right\vert
<\varepsilon $: 
\begin{equation*}
\text{\qquad }\frac{\partial r(\zeta )\wedge \overline{\partial }r(\zeta
)\wedge \partial \overline{\partial }r(\zeta )}{\Phi _{K}(\zeta ,z)}=\omega
_{n}\wedge \overline{\omega }_{n}\wedge \left[ \sum_{j=1}^{n-1}{}A_{j}d\zeta
_{j}\wedge \overline{d\zeta _{j}}+\sum_{j,l}{}B_{jl}\ d\zeta _{j}\wedge 
\overline{d\zeta _{l}}\right] ,
\end{equation*}%
where%
\begin{equation*}
\left\vert A_{j}(\zeta ,z)\right\vert \lesssim \frac{1}{\left\vert \func{Im}%
F^{(r)}(\zeta ,z)\right\vert +\left\vert r(\zeta )\right\vert +\left\vert
r(z)\right\vert +\left\vert \zeta _{j}-z_{j}\right\vert ^{2}+\frac{K}{2}%
\left\vert \zeta -z\right\vert ^{3}~}
\end{equation*}%
and%
\begin{equation*}
\left\vert B_{jl}(\zeta ,z)\right\vert \lesssim \frac{1}{\left\vert \func{Im}%
F^{(r)}(\zeta ,z)\right\vert +\left\vert r(\zeta )\right\vert +\left\vert
r(z)\right\vert +\frac{K}{2}\left\vert \zeta -z\right\vert ^{2}}\text{.}
\end{equation*}
\end{corollary}

\textit{Proof. }From the preceding propositions one obtains%
\begin{equation*}
\frac{\partial r\wedge \overline{\partial }r\wedge \partial \overline{%
\partial }r(\zeta )}{\Phi _{K}(\zeta ,z)}=\partial r\wedge \overline{%
\partial }r\wedge \left[ \frac{1/2\sum_{j=1}^{n-1}\lambda _{j}(z)d\zeta
_{j}\wedge \overline{d\zeta _{j}}}{\Phi _{K}(\zeta ,z)}+\frac{\Omega _{1}}{%
\Phi _{K}(\zeta ,z)}\right] .
\end{equation*}%
By using the estimate for $\left\vert \Phi _{K}\right\vert $ in Theorem 3
and after cancelling $\left\vert \zeta -z\right\vert $, one readily sees
that the coefficients $B_{jl}$ of the form $\Omega _{1}/\Phi _{K}$ satisfy
the required estimate. For the leading terms, estimate 
\begin{equation*}
\left\vert \frac{\lambda _{j}(z)}{\Phi _{K}(\zeta ,z)}\right\vert \lesssim 
\frac{\lambda _{j}(z)}{\left\vert \func{Im}F^{(r)}(\zeta ,z)\right\vert
+\left\vert r(\zeta )\right\vert +\left\vert r(z)\right\vert +\lambda
_{j}(z)\left\vert \zeta _{j}-z_{j}\right\vert ^{2}+\frac{K}{2}\left\vert
\zeta -z\right\vert ^{3}}.
\end{equation*}%
Observe that there exists a constant $C>0$ independent of $z$, such that $%
0\leq $ $\lambda _{j}(z)\leq C$ for $j=1,...,n-1$. Thus, if $\lambda
_{j}(z)>0$, one has $1/\lambda _{j}(z)\geq 1/C>0$, and one may cancel the
factor $\lambda _{j}(z)$ from numerator and denominator to obtain (with a
new constant in $\lesssim $)%
\begin{equation*}
\left\vert \frac{\lambda _{j}(z)}{\Phi _{K}(\zeta ,z)}\right\vert \lesssim 
\frac{1}{\left\vert \func{Im}F^{(r)}(\zeta ,z)\right\vert +\left\vert
r(\zeta )\right\vert +\left\vert r(z)\right\vert +\left\vert \zeta
_{j}-z_{j}\right\vert ^{2}+\frac{K}{2}\left\vert \zeta -z\right\vert ^{3}}.
\end{equation*}%
Trivially this estimate holds also $\lambda _{j}(z)=0$. $\blacksquare $

\textbf{Remark.}\textit{\ }The special $z$-diagonalizing coordinates in
dependence of $z$ introduced above are used only for the estimations in the
preceding propositions and corollary, and for relevant estimates in section
5. Unless something different is explicitly mentioned, it is assumed that in
all other places the coordinates are the fixed standard coordinates of $%
\mathbb{C}^{n}$ introduced at the beginning. In particular, the function $%
\Phi _{K}$ is defined with respect to these standard coordinates of $\mathbb{%
C}^{n}$, which will continue to be used in the following sections.

\section{Factorization and the generating form}

Following standard procedure, in order to build a kernel generating form on $%
bD\times (D\cap U)$ from the support function $\Phi _{K}$, one needs a
decomposition%
\begin{equation*}
\Phi _{K}(\zeta ,z)=\sum_{j=1}^{n}g_{j}(\zeta ,z)(\zeta _{j}-z_{j})\text{.}
\end{equation*}%
Note that for $\zeta \in bD$ one has $\Phi _{K}(\zeta ,z)=\Phi _{K}(\zeta
,z)+r(\zeta )=F^{(r)}(\zeta ,z)+K\left\vert \zeta -z\right\vert ^{3}$.
Clearly there is such a decomposition for the Levi polynomial $F^{(r)}$,
with the corresponding coefficients holomorphic in $z.$ Since 
\begin{equation*}
\left\vert \zeta -z\right\vert ^{3}=\sum_{j=1}^{n}\left\vert \zeta
-z\right\vert \overline{(\zeta _{j}-z_{j})}(\zeta _{j}-z_{j}),
\end{equation*}%
there is a decomposition for $\Phi _{K}+r(\zeta )$ with $g_{j}$ equal to a
sum of a holomorphic term and a term of type $\mathcal{E}_{2}^{\#}$. This
implies the following lemma.

\begin{lemma}
There is a decomposition 
\begin{equation*}
\Phi _{K}(\zeta ,z)+r(\zeta )=\sum_{j=1}^{n}g_{j}(\zeta ,z)(\zeta _{j}-z_{j})%
\text{ for }\zeta \in U\text{,}
\end{equation*}%
where $\overline{\partial _{z}}g_{j}$ is of type $\mathcal{E}_{1}^{\#}$ for $%
j=1,...,n$.
\end{lemma}

We now define the $(1,0)$ form $g=\sum_{j=1}^{n}g_{j}d\zeta _{j}$ and set 
\begin{equation*}
W^{K}(\zeta ,z)=\frac{g}{\Phi _{K}(\zeta ,z)}
\end{equation*}%
for $(\zeta ,z)\in \left( \overline{D}\cap U\right) \times \left( \overline{D%
}\cap U\right) $ with $0<\left\vert \zeta -z\right\vert <\varepsilon $. Note
that the $(1,0)$ form $W^{K}(\zeta ,z)$ has smooth coefficients for $\zeta
\neq z$ and that it satisfies 
\begin{equation*}
<W^{K}(\zeta ,z),\zeta -z>=\sum_{j=1}^{n}(g_{j}/\Phi _{K})(\zeta
_{j}-z_{j})=1
\end{equation*}%
for $\zeta \in bD$ and $0<\left\vert \zeta -z\right\vert <\varepsilon $; it
therefore is a (local) generating form in the terminology of Ra86.

In order to globalize $\Phi _{K}$ and $W^{K}$ in $z$, we patch with the
corresponding terms from the Bochner-Martinelli kernel, as was done in Ra86
in case of strictly pseudoconvex domains. We choose a $C^{\infty }$ function 
$\chi (t)$ such that $0\leq \chi (t)\leq 1,$ $\chi (t)=1$ for $t\leq
\varepsilon /2$, and $\chi (t)=0$ for $t\geq 3/4\varepsilon $. Define $%
S(\zeta ,z)$ on $U\times \overline{D}$ by 
\begin{equation*}
S(\zeta ,z)=\chi (\left\vert \zeta -z\right\vert )\Phi _{K}(\zeta
,z)+[1-\chi (\left\vert \zeta -z\right\vert )]\left\vert \zeta -z\right\vert
^{2}\text{.}
\end{equation*}%
We need to ensure that $S$ does not have any new zeroes.

\begin{lemma}
$S(\zeta ,z)\neq 0$ for $\zeta \in \overline{D}\cap U$ with $\left\vert
r(\zeta )\right\vert \leq \varepsilon $ and all $z\in \overline{D}$ with $%
z\neq \zeta .$
\end{lemma}

\textit{Proof.} Since 
\begin{eqnarray*}
\left\vert S\right\vert &\succeq &\left\vert \func{Re}S\right\vert
+\left\vert \func{Im}S\right\vert \geq \func{Re}S+\left\vert \func{Im}%
S\right\vert \geq \\
&\geq &\chi \func{Re}\Phi _{K}(\zeta ,z)+(1-\chi )\left\vert \zeta
-z\right\vert ^{2}\text{,}
\end{eqnarray*}%
and since $\chi (\left\vert \zeta -z\right\vert )\equiv 0$ for $\left\vert
\zeta -z\right\vert \geq 3/4\varepsilon $, we can estimate the first term by
utilizing the estimate for $\func{Re}\Phi _{K}(\zeta ,z)$ from section 3 to
obtain 
\begin{equation*}
\left\vert S\right\vert \succeq \chi {\Large [}\left\vert r(\zeta
)\right\vert /2+\left\vert r(z)\right\vert +\mathcal{L}(r,\zeta ;\pi _{\zeta
}^{t}(\zeta -z))+K\left\vert \zeta -z\right\vert ^{3}{\Large ]}-\chi
\left\vert \Phi _{K}\right\vert /2+(1-\chi )\left\vert \zeta -z\right\vert
^{2}\text{.}
\end{equation*}%
By estimating $\chi \left\vert \Phi _{K}\right\vert \leq \left\vert
S\right\vert +(1-\chi )\left\vert \zeta -z\right\vert ^{2}$ it then follows
that 
\begin{equation*}
\left\vert S\right\vert \succeq \chi K\left\vert \zeta -z\right\vert ^{3}-1/2%
{\Large [}\left\vert S\right\vert +(1-\chi )\left\vert \zeta -z\right\vert
^{2}{\Large ]}+(1-\chi )\left\vert \zeta -z\right\vert ^{2}.
\end{equation*}%
Therefore%
\begin{equation*}
\left\vert S\right\vert \succeq \chi K\left\vert \zeta -z\right\vert
^{3}+1/2(1-\chi )\left\vert \zeta -z\right\vert ^{2}\succeq \left\vert \zeta
-z\right\vert ^{3},
\end{equation*}%
and the lemma is proved.

For $(\zeta ,z)\in bD\times \overline{D}$ one has the decomposition%
\begin{equation*}
S(\zeta ,z)=\sum_{j=1}^{n}s_{j}(\zeta ,z)(\zeta _{j}-z_{j})\text{, }
\end{equation*}%
with $s_{j}(\zeta ,z)=\chi (\left\vert \zeta -z\right\vert )g_{j}(\zeta
,z)+(1-\chi )(\overline{\zeta _{j}-z_{j}}).$ Finally, by introducing the $%
(1,0)$ form $s=\sum_{j}s_{j}(\zeta ,z)~d\zeta _{j}$, one obtains the global
generating form%
\begin{equation*}
W^{S}(\zeta ,z)=\frac{s(\zeta ,z)}{S(\zeta ,z)}
\end{equation*}%
on $bD\times \overline{D}-\{(\zeta ,\zeta ):\zeta \in bD\}$. Note that for $%
\zeta \in bD$ and $0<\left\vert \zeta -z\right\vert \leq \varepsilon /2$ one
has $W^{S}(\zeta ,z)=g/\Phi _{K}$.

We now introduce the corresponding Cauchy-Fantappi\'{e} kernel (of order $0)$
\begin{equation*}
\Omega _{0}(W^{S})=(2\pi i)^{-n}W^{S}\wedge (\overline{\partial _{\zeta }}%
W^{S})^{n-1}\text{ on }bD\times \overline{D}-\{(\zeta ,\zeta ):\zeta \in bD\}%
\text{.}
\end{equation*}%
We recall the following standard properties of Cauchy-Fantappi\'{e} kernels,
that is, for their pull-backs to $bD$:

\begin{equation*}
\text{a)}\overline{\text{ }\partial _{\zeta }}\Omega _{0}(W^{S})=0,\text{ 
\hspace{0.5in} b) }\Omega _{0}(W^{S})=(2\pi i)^{-n}\frac{s\wedge (\overline{%
\partial _{\zeta }}s)^{n-1}}{S^{n}}\text{,}
\end{equation*}%
and 
\begin{equation*}
\text{c) }f(z)=\int_{bD}f(\zeta )~\Omega _{0}(W^{S})~\text{for }f\in 
\mathcal{O}(D)\cap C(\overline{D})\text{ and }z\in D\text{.}
\end{equation*}%
(See Ra86, for example.)

\section{Some regularity properties}

As we pointed out at the beginning, for a given weakly pseudoconvex domain
it is in general not possible to find an explicit Cauchy-Fantappi\'{e}
kernel which is holomorphic near the singularity $z=\zeta \in bD$. The
kernel $\Omega _{0}(W^{S})$ we constructed in the preceding section is---in
some way---optimal in this general setting. On the one hand, just like the
universal Bochner-Martinelli kernel, if $z\in bD,$ the singularity of $%
\Omega _{0}(W^{S})$ at $\zeta =z$ is not integrable over $bD$, but it is
right at the border line: as we will show, for any $\alpha >0,$ $\left\vert
\zeta -z\right\vert ^{\alpha }\Omega _{0}(W^{S})$ is indeed integrable over $%
bD$. This implies that the operator $T^{S}:$ $L^{1}(bD)\rightarrow C(D)$
defined by 
\begin{equation*}
T^{S}(f)=\int_{bD}f(\zeta )~\Omega _{0}(W^{S})
\end{equation*}
preserves a variety of function spaces, at least if one allows for an
arbitrarily small loss in regularity. Such results, however, require
detailed new proofs, as the singularity of the operator $T^{S}$ does not
appear to fit into any of the classical theories. Much more significant is
the fact that $\Omega _{0}(W^{S})$ reflects the \emph{complex} geometry of
the boundary $bD$, and hence should be much more useful in complex analysis
than the Bochner-Martinelli kernel. For example, while $T^{S}(f)$ is not
holomorphic in general, one has the following result, which makes it
explicit that $T^{S}(f)$ enjoys some special complex analytic properties.

\begin{theorem}
\bigskip For any $\delta <2/3$ there exists a constant $C_{\delta }$ such
that for all functions $f$ continuous on $bD$ one has $T^{S}(f)\in C^{\infty
}(D)$ and 
\begin{equation*}
\left\vert \overline{\partial _{z}}T^{S}(f)(z)\right\vert \leq C_{\delta
}\left\vert f\right\vert _{0}dist(z,bD)^{\delta -1}\text{ for }z\in D.%
\footnote{$\left\vert f\right\vert _{0}$ denotes the supremum norm of $f$
over $bD$.}
\end{equation*}
\end{theorem}

In contrast, the Bochner-Martinelli kernel satisfies the analogous estimate
only for $\delta =0$; it does not give preference to partial derivatives
with respect to the \emph{conjugate} variables $\overline{z_{j}}$.

We note that given the structure of the coefficients of $\Omega _{0}(W^{S})$
it follows by standard arguments that $T^{S}(f)\in C^{\infty }(D)$ for any $%
f $ that is integrable over $bD.$

For the proof of the desired estimate we may differentiate with respect to $%
\overline{z_{j}}$ under the integral sign, and we shall utilize variations
of the estimates for $\overline{\partial }\partial r/\Phi _{K}$ that were
established in section 3. Since $\left\vert S\right\vert \geq \gamma >0$ for 
$\left\vert \zeta -z\right\vert \geq \varepsilon /2$, for the purposes of
estimations it is enough to consider the case where $\zeta \in bD$ and $%
\left\vert \zeta -z\right\vert \leq \varepsilon /2,$ so that 
\begin{equation*}
s_{j}=g_{j}=\partial r/\partial \zeta _{j}-1/2\sum_{k}\partial
^{2}r/\partial \zeta _{j}\partial \zeta _{k}(\zeta _{k}-z_{k})+\mathcal{E}%
_{2}^{\#}\text{.}
\end{equation*}%
By making use of the specific form of the second term of $s_{j}=g_{j}$, it
follows that 
\begin{eqnarray*}
s &=&\partial r(\zeta )+\mathcal{E}_{1}^{\#}\text{, and} \\
\overline{\partial }_{\zeta }s &=&\overline{\partial }\partial r(\zeta )+%
\mathcal{E}_{1}^{\#}\text{.}
\end{eqnarray*}%
We now fix $z\in U\cap \overline{D}$ and introduce the frame $%
L_{1},...,L_{n} $ and the $z$-diagonalizing coordinates for $\zeta $, as in
the proof of Proposition 4. The proof of Proposition 5 shows that%
\begin{equation*}
\overline{\partial _{\zeta }}s=\left[ \mathcal{L(}z)~\mathcal{+E}_{1}^{\#}%
\right] +N\text{,}
\end{equation*}%
where $\mathcal{L}(z)\mathcal{=}1/2\sum_{j=1}^{n-1}\lambda _{j}(z)d\zeta
_{j}\wedge \overline{d\zeta _{j}}$, and $N$ represents a "normal" $2$-form
of the type $\Lambda _{1}\wedge \omega _{n}+\Lambda _{2}\wedge \overline{%
\omega _{n}}${\large ,} where $\Lambda _{1}$ and $\Lambda _{2}$ are suitable 
$1$-forms which may change from formula to formula. Note that the pull-back
of $N\wedge N$ to $bD$ is zero; we shall therefore ignore such terms in the
following.

\begin{lemma}
If $\iota _{bD}:$ $bD\rightarrow U$ is the inclusion, for any $1\leq q\leq
n-1$ the form $\iota _{bD}^{\ast }(\partial s\wedge (\overline{\partial
_{\zeta }}s)^{q}$ is a linear combination with $\mathcal{E}_{0}^{\#}$
coefficients of $(q+1,q)$-forms which contain factors of the type 
\begin{equation*}
\sum_{\left\vert J\right\vert =l}{}^{^{\prime }}\lambda _{J}(z)(d\zeta
\wedge \overline{d\zeta })^{J}\wedge \lbrack \mathcal{E}_{1}^{\#}]^{q-l}%
\text{ with }0\leq l\leq q\text{,}
\end{equation*}%
where the summation is over strictly increasing $l$-tuples $%
J=[j_{1},...,j_{l}]\subset $\newline
$[1,2...,n-1]$, 
\begin{equation*}
\lambda _{J}(z)(d\zeta \wedge \overline{d\zeta })^{J}=\prod\nolimits_{\nu
=1}^{l}\lambda _{j_{\nu }}(z)d\zeta _{j_{v}}\wedge \overline{d\zeta _{j_{v}}}%
\text{,}
\end{equation*}%
and $[\mathcal{E}_{1}^{\#}]^{q-l}$ is a form of degree $2(q-l)$ whose
coefficients are products of $q-l$ factors of type $\mathcal{E}_{1}^{\#}$.
\end{lemma}

\emph{Proof.} By the above representations for $s$ and $\overline{\partial
_{\zeta }}s,$ one obtains%
\begin{equation*}
s\wedge (\overline{\partial _{\zeta }}s)^{q}=N\wedge (\overline{\partial
_{\zeta }}s)^{q}+\mathcal{E}_{1}^{\#}\wedge (\overline{\partial _{\zeta }}%
s)^{q}
\end{equation*}%
and 
\begin{eqnarray*}
(\overline{\partial _{\zeta }}s)^{q} &=&\left[ \mathcal{L(}z)~\mathcal{+E}%
_{1}^{\#}\right] ^{q}+\left[ \mathcal{L(}z)~\mathcal{+E}_{1}^{\#}\right]
^{q-1}\wedge N= \\
&=&\sum_{l=0}^{q}[\mathcal{L(}z)]^{l}\wedge ~[\mathcal{E}_{1}^{\#}]^{q-l}+%
\sum_{l=0}^{q-1}[\mathcal{L(}z)]^{l}\wedge ~[\mathcal{E}_{1}^{\#}]^{q-1-l}%
\wedge N\text{,}
\end{eqnarray*}%
where we have ignored terms with $N\wedge N$. It follows that%
\begin{multline*}
s\wedge (\overline{\partial _{\zeta }}s)^{q}=(N+\mathcal{E}_{1}^{\#})\wedge
\sum_{l=0}^{q}[\mathcal{L(}z)]^{l}\wedge ~[\mathcal{E}_{1}^{\#}]^{q-l}+ \\
+\mathcal{E}_{1}^{\#}\wedge \sum_{l=0}^{q-1}[\mathcal{L(}z)]^{l}\wedge ~[%
\mathcal{E}_{1}^{\#}]^{q-1-l}\wedge N=
\end{multline*}
\begin{equation*}
=N\wedge \sum_{l=0}^{q}[\mathcal{L(}z)]^{l}\wedge ~[\mathcal{E}%
_{1}^{\#}]^{q-l}+.....\text{,}
\end{equation*}%
where we have retained just the leading terms, i.e., those with the lowest
order of vanishing. The lemma then follows by expanding $[\mathcal{L(}%
z)]^{l}={\LARGE [}1/2\sum_{j=1}^{n-1}\lambda _{j}(z)d\zeta _{j}\wedge 
\overline{d\zeta _{j}}{\Large ]}^{l}$ and absorbing all numerical factors in
the $N$ and/or $\mathcal{E}_{1}^{\#}$ terms.

Before coming to the proof of theorem 9, we use the preceding results to
prove the following statement that we had anticipated earlier.

\begin{proposition}
For any $\alpha >0$ there exists a constant $C_{\alpha }$ such that for any $%
z\in \overline{D}$ with $dist(z,bD)<\varepsilon /2$ one has%
\begin{equation*}
\int_{bD}\left\vert \zeta -z\right\vert ^{\alpha }\left\vert \Omega
_{0}(W^{S})(\zeta ,z)\right\vert \leq C_{a}\text{.}
\end{equation*}
\end{proposition}

\textit{Proof. }We fix $z$ as required and introduce the $z$-diagonalizing
coordinates on a ball $B(z,\delta )$, with $\delta \leq \varepsilon /2$.%
\textit{\ }It is clearly enough to prove the estimate locally, i.e., for $%
\zeta \in bD\cap B(z,\delta )$. From Lemma 10, with $q=n-1$, one sees that
the pull-back to $bD$ of $\left\vert \zeta -z\right\vert ^{\alpha }\Omega
_{0}(W^{S})=$ $\left\vert \zeta -z\right\vert ^{\alpha }c_{n}s\wedge (%
\overline{\partial _{\zeta }}s)^{n-1}/\Phi _{K}^{n}$ decomposes into a
linear combination of terms

\begin{equation*}
\left\vert \zeta -z\right\vert ^{\alpha }\frac{\omega _{n}\wedge
\sum_{\left\vert J\right\vert =l}{}^{^{\prime }}\lambda _{J}(z)(d\zeta
\wedge \overline{d\zeta })^{J}\wedge \lbrack \mathcal{E}_{1}^{\#}]^{n-1-l}}{%
\Phi _{K}^{n}}
\end{equation*}%
for $l=0,...,n-1$, and of other terms of this type with higher order of
vanishing in the numerator. In order to integrate over $bD\cap B(z,\delta )$
we introduce a $C^{k-2}$ coordinate system $(\zeta _{1}-z_{1},...,\zeta
_{n-1}-z_{n-1},\func{Im}F^{r}(\zeta ,z),r\left( \zeta \right) )$ on $%
B(z,\delta )$, where $\delta $ may have to be chosen smaller, but can in any
case be fixed independently of the point $z$. After renumbering, and since $%
r(\zeta )=0$ on $bD$, one is left with estimating%
\begin{eqnarray*}
I_{l}(z) &=&\int_{bD\cap B(z,\delta )}\left\vert \zeta -z\right\vert
^{\alpha }\frac{\lambda _{1}(z)...\lambda _{l}(z)\left\vert \zeta
-z\right\vert ^{n-1-l}}{\left\vert \Phi _{K}\right\vert ^{n}}\cdot \\
&&\cdot \left\vert d\zeta _{1}\wedge d\overline{\zeta _{1}}\wedge ...\wedge
d\zeta _{n-1}\wedge d\overline{\zeta _{n-1}}\wedge d\func{Im}F^{r}(\zeta
,z)\right\vert
\end{eqnarray*}%
for any $l=0,...,n-1$. Choose $\gamma >0$ so that $\gamma ^{\#}=\alpha
-3\gamma (l+1)>0$, and set $\zeta ^{^{\prime }}=(\zeta _{l+1},...,\zeta
_{n-1})\in \mathbb{C}^{n-1-l}$. Since $\left\vert \Phi _{K}\right\vert
\succeq \left\vert \zeta -z\right\vert ^{3}$ one has $\left\vert \zeta
-z\right\vert ^{\alpha }/\left\vert \Phi _{K}\right\vert ^{\gamma
(l+1)}\preceq \left\vert \zeta -z\right\vert ^{\gamma ^{\#}}$, and it
follows from Proposition 4 that 
\begin{equation*}
\lambda _{j}(z)/\left\vert \Phi _{K}\right\vert ^{1-\gamma }\lesssim
1/[\left\vert r(z)\right\vert +\left\vert \zeta _{j}-z_{j}\right\vert
^{2}]^{1-\gamma }.
\end{equation*}
By using these estimates one obtains 
\begin{eqnarray*}
I_{l}(z) &\preceq &\int_{bD\cap B(z,\delta )}\prod\nolimits_{j=1}^{l}\frac{%
\left\vert d\zeta _{j}\wedge d\overline{\zeta _{j}}\right\vert }{[\left\vert
r(z)\right\vert +\left\vert \zeta _{j}-z_{j}\right\vert ^{2}]^{1-\gamma }}%
\cdot \frac{d\func{Im}F^{r}}{\left\vert \func{Im}F^{r}\right\vert ^{1-\gamma
}} \\
&&\cdot \frac{\left\vert \zeta -z\right\vert ^{\gamma ^{\#}+n-1-l}}{%
\left\vert \Phi _{K}\right\vert ^{(n-1-l)}}dV(\zeta ^{\prime })\text{.}
\end{eqnarray*}%
Note that in case $l=n-1$ the last factor is unifomly bounded, while if $%
l<n-1$, one has 
\begin{equation*}
\frac{\left\vert \zeta -z\right\vert ^{\gamma ^{\#}+n-1-l}}{\left\vert \Phi
_{K}\right\vert ^{(n-1-l)}}\preceq \frac{1}{\left\vert \zeta ^{\prime
}-z^{\prime }\right\vert ^{2(n-1-l)-\gamma ^{\#}}}\text{ ,}
\end{equation*}%
which is integrable in $\zeta ^{^{\prime }}.$ Consequently all factors in
the resulting iterated integral for $I_{l}(z)$ are bounded independently of $%
z$.$\blacksquare $

As an application, one obtains the following result for the space $\Lambda
_{\alpha }(bD)$ of functions which are H\"{o}lder continuous of order $%
\alpha $, with norm $\left\vert \cdot \right\vert _{\alpha }$.

\begin{corollary}
For any $\alpha >0$ the integral operator $T^{S}:L^{1}(bD)\rightarrow
C^{\infty }(D)$ satisfies%
\begin{equation*}
\left\vert T^{S}(f)\right\vert _{0}\leq C_{\alpha }\left\vert f\right\vert
_{\alpha }
\end{equation*}%
for all $f\in \Lambda _{\alpha }(bD)$.
\end{corollary}

More generally, it is quite likely that for any $0<\alpha ^{\prime }<\alpha
<1$ the operator $T^{S}$ is bounded from $\Lambda _{\alpha }(bD)$ to $%
\Lambda _{\alpha ^{\prime }}(D).$ Such variations will be discussed at some
other time.

\textit{Proof of the theorem. }Since $\left\vert r\left( z\right)
\right\vert \approx dist(z,bD)$, it is enough to prove the estimate in the
theorem with $\left\vert r\left( z\right) \right\vert $ instead of $%
dist(z,bD)$. We apply $\overline{\partial _{z}}$ to $T^{S}(f)$ under the
integral sign. Note that 
\begin{eqnarray*}
\overline{\partial _{z}}\Omega _{0}(W^{S}) &=&c_{n}\frac{\overline{\partial
_{z}}s\wedge (\overline{\partial _{\zeta }}s)^{n-1}+(n-1)~s\wedge \overline{%
\partial _{z}}\overline{\partial _{\zeta }}s\wedge (\overline{\partial
_{\zeta }}s)^{n-2}}{\Phi _{K}^{n}}+ \\
&&-c_{n}n\frac{s\wedge (\overline{\partial _{\zeta }}s)^{n-1}\wedge 
\overline{\partial _{z}}\Phi _{K}}{\Phi _{K}^{n+1}}.
\end{eqnarray*}%
As in the proof of the proposition, we fix $z$ and choose the frame and
unitary coordinate change adapted to that point $z.$ By analyzing the terms
above as before, and estimating $\overline{\partial _{z}}\Phi _{K}=\mathcal{E%
}_{2}^{\#}$ by proposition 2, one sees that the critical integrals that one
needs to estimate are those with factors of type

\begin{equation*}
\text{(I) }\frac{\omega _{n}\wedge \sum_{\left\vert J\right\vert
=l}{}^{^{\prime }}\lambda _{J}(z)(d\zeta \wedge d\overline{\zeta }%
)^{J}\wedge \lbrack \mathcal{E}_{1}^{\#}]^{n-2-l}}{\Phi _{K}^{n}}\text{ for }%
0\leq l\leq n-2
\end{equation*}%
and 
\begin{equation*}
\text{(II) }\frac{\omega _{n}\wedge \sum_{\left\vert J\right\vert
=l}{}^{^{\prime }}\lambda _{J}(z)(d\zeta \wedge d\overline{\zeta }%
)^{J}\wedge \lbrack \mathcal{E}_{1}^{\#}]^{n-1-l}}{\Phi _{K}^{n}}\frac{%
\mathcal{E}_{2}^{\#}}{\Phi _{K}}\text{ for }0\leq l\leq n-1.
\end{equation*}%
Note that $\left\vert \mathcal{E}_{1}^{\#}\mathcal{E}_{2}^{\#}/\Phi
_{K}\right\vert \ $is bounded by a constant. Therefore, if $l\leq n-2,$ the
terms in (II) are estimated by this constant multiplied with terms of type
(I). We thus need to estimate terms of type (I), and those of type (II) with 
$l=n-1.$ Consider (I) first. For a fixed $l\leq n-2$---after renumbering
indices---it is enough to consider 
\begin{gather*}
I_{l}(z)=\int_{bD\cap B(z,\delta /2)}\frac{\lambda _{1}(z)...\lambda
_{l}(z)\left\vert \zeta -z\right\vert ^{n-2-l}}{\left\vert \Phi
_{K}\right\vert ^{n}}~ \\
\hspace{0.7in}\hspace{0.7in}\cdot \left\vert d\zeta _{1}\wedge d\overline{%
\zeta _{1}}\wedge ...\wedge d\zeta _{n-1}\wedge d\overline{\zeta _{n-1}}%
\wedge d\func{Im}F^{r}(\zeta ,z)\right\vert \text{ }
\end{gather*}%
for $z\in D\cap U$. Fix $0<\delta <2/3$, and choose $\gamma >0$ so small
that $3(l+1)\gamma <2-3\delta $. We factor $\left\vert \Phi _{K}\right\vert
^{l+1}=\left\vert \Phi _{K}\right\vert ^{(l+1)(1-\gamma )}\left\vert \Phi
_{K}\right\vert ^{(l+1)\gamma }$ as before and consider the remaining factor 
\begin{equation*}
J_{l}(z)=\frac{\left\vert \zeta -z\right\vert ^{n-2-l}}{\left\vert \Phi
_{K}\right\vert ^{n-1-l+(l+1)\gamma }}\text{.}
\end{equation*}%
Let us use the estimate $\left\vert \Phi _{K}\right\vert \succeq \left\vert
\Phi _{K}\right\vert ^{\delta }\left\vert r\left( z\right) \right\vert
^{1-\delta }\succeq $ $\left\vert \zeta -z\right\vert ^{3\delta }\left\vert
r\left( z\right) \right\vert ^{1-\delta }$ for one of the factors in the
denominator, and the estimate $\left\vert \Phi _{K}\right\vert \succeq $ $%
\left\vert \zeta -z\right\vert ^{3}$ for the other factors. It follows that 
\begin{eqnarray*}
J_{l}(z) &\lesssim &\left\vert r\left( z\right) \right\vert ^{\delta -1}%
\frac{\left\vert \zeta -z\right\vert ^{n-2-l}}{\left\vert \zeta
-z\right\vert ^{3[n-2-l+\delta +(l+1)\gamma ]}}\preceq \\
&\preceq &\left\vert r\left( z\right) \right\vert ^{\delta -1}\frac{1}{%
\left\vert \zeta ^{\prime }-z^{\prime }\right\vert ^{2[n-2-l]+3[\delta
+(l+1)\gamma ]}}\text{,}
\end{eqnarray*}%
where $\zeta ^{^{\prime }}=(\zeta _{l+1},...,\zeta _{n-1})\in \mathbb{C}%
^{n-1-l}$. Since by the choice of $\gamma $ one has $3[\delta +(l+1)\gamma
]<2$, this last expression is integrable in $\zeta ^{\prime }\in B(z^{\prime
},\delta /2).$ It follows that the integral $I_{l}(z)$ is bounded by $%
\left\vert r\left( z\right) \right\vert ^{\delta -1}$ multiplied with a
uniformly bounded iterated integral. Finally, we must consider the term (II)
for $l=n-1.$ Choose $\gamma >0$ so that $3n\gamma \leq 2-3\delta $, and
factor $\left\vert \Phi _{K}\right\vert ^{n}=\left\vert \Phi _{K}\right\vert
^{n(1-\gamma )}\left\vert \Phi _{K}\right\vert ^{n\gamma }$. By factoring
the integrand as before, and by estimating the remaining factor 
\begin{equation*}
\frac{\left\vert \zeta -z\right\vert ^{2}}{\left\vert \Phi _{K}\right\vert
^{1+n\gamma }}\lesssim \frac{\left\vert \zeta -z\right\vert ^{2}}{\left\vert
\Phi _{K}\right\vert ^{\delta +n\gamma }\left\vert r\left( z\right)
\right\vert ^{1-\delta }}\lesssim \frac{\left\vert \zeta -z\right\vert ^{2}}{%
\left\vert \zeta -z\right\vert ^{3(\delta +n\gamma )}\left\vert r\left(
z\right) \right\vert ^{1-\delta }}\preceq \left\vert r\left( z\right)
\right\vert ^{\delta -1},
\end{equation*}%
the desired estimate follows.$\blacksquare $

\section{Concluding remarks\textbf{\ }}

The estimates we proved in the preceding section show that the support
function $\Phi _{K}$, the generating form $W^{S}=s/S,$ and the associated
Cauchy-Fantappi\'{e} kernel $\Omega _{0}(W^{S})$ provide a \textit{partial}
replacement for the missing holomorphic Cauchy kernel on arbitrary weakly
pseudoconvex domains. Similarly, the corresponding higher order forms $%
\Omega _{q}(W^{S})$, $0\leq q\leq n-1$, have proved useful in applications
to pointwise estimates in the $\overline{\partial }-$Neumann theory (see
Ra11), and as indicated there, might provide critical basic ingredients to
prove deep new results in case the domain is assumed to be of \emph{finite}
type. Such applications, along with others, will be the subject of
forthcoming articles.\bigskip

\section*{References}

.

Ca87\qquad Catlin, D.: Subelliptic estimates for the $\overline{\partial }-$%
Neumann problem on pseudoconvex domains. \textit{Ann. Math. }\textbf{126 }%
(1987), 131- 191.

Cu97\qquad Cumenge, A.: Estim\'{e}es Lipschitz optimales dans les convexes
de type fini. \textit{C. R. Acad. Sci. Paris} \textbf{325} (1997), 1077-1080.

Da82\qquad D'Angelo, J. F.: Real hypersurfaces, orders of contact, and
applications. \textit{Ann. Math.} \textbf{115} (1982), 615-637.

DF78\qquad Diederich, K., and Fornaess, J. E.: Pseudoconvex domains with
real-analytic boundary. \textit{Ann. Math.} \textbf{107} (1978), 371-384.

DF99\qquad Diederich, K., and Fornaess, J. E.: Support functions for convex
domains of finite type. \textit{Math. Z.} \textbf{230} (1999), 145-164.

DFF99 Diederich, K., Fischer, B., and Fornaess, J. E.: H\"{o}lder estimates
on convex domains of finite type. \textit{Math. Z.} \textbf{232} (1999), 43
- 61.

FoKo72\qquad Folland, G., and Kohn, J.J.: \textit{The Neumann Problem for
the Cauchy-Riemann Complex.} Princeton Univ. Press, 1972.

Ko72\qquad Kohn, J.J.: Boundary behavior of $\overline{\partial }$ on weakly
pseudoconvex manifolds of dimension two. \textit{J. Diff. Geometry} \textbf{6%
} (1972), 523 -542.

Ko79 \qquad Kohn, J.J.: Subelllipticity of the $\overline{\partial }-$%
Neumann Problem on Pseudoconvex Domains: Sufficient Conditions. \textit{Acta
math.} \textbf{142} (1979), 79-122.

KoNi72\qquad Kohn, J.J., and Nirenberg, L.: A pseudoconvex domain not
admitting a holomorphic support function. \textit{Math. Ann. }\textbf{201 }%
(1973), 265 - 268.

Ra78\qquad Range, R.M.: On H\"{o}lder estimates for $\overline{\partial }u=f$
on weakly pseudoconvex domains. \textit{Proc. Int. Conf. Cortona 1976-77}.
Sc. Norm. Sup. Pisa (1978), 247-267.

Ra81\qquad A remark on bounded strictly plurisubharmonic exhaustion
functions. \textit{Proc. Amer. Math. Soc. }\textbf{81 }(1981), 220 - 222.

Ra86\qquad Range, R. M.: \textit{Holomorphic Functions and Integral
Representations in Several Complex Variables.} Springer Verlag New York
1986, corrected 2nd. printing 1998.

Ra90\qquad Range, R. M.: Integral kernels and H\"{o}lder estimates for $%
\overline{\partial }$ on pseudoconvex domains of finite type in $\mathbb{C}%
^{2}$. \textit{Math. Ann.} \textbf{288} (1990), 63-74.

Ra11\qquad Range, R. M.: A pointwise basic estimate and H\"{o}lder
multipliers for the $\overline{\partial }-$Neumann problem on pseudoconvex
domains. \textit{Preliminary Report.} \textit{arXiv:1106.3132 (June 2011)}

Si10 \qquad Siu, Y.-T.: Effective termination of Kohn's algorithm for
subelliptic multipliers. \textit{Pure Appl. Math. Quarterly} \textbf{6}
(2010), no. 4, Special Issue: In honor of Joseph J. Kohn.\vspace{0.5in}

R. Michael Range

Department of Mathematics

State University of New York at Albany

range@albany.edu

\end{document}